\documentclass[a4paper,american,reqno,oneside]{amsart}

\usepackage{amssymb,amsmath,amsfonts,amsthm}
\usepackage{booktabs}
\usepackage{algorithm}
\usepackage[noend]{algorithmic}
\usepackage{graphicx}
\usepackage{pgfplots}
\usetikzlibrary{pgfplots.groupplots}
\usepackage[binary-units=true]{siunitx}
\usepackage[style=authoryear-comp,
            giveninits=true,
            dashed=false,
            maxbibnames=100,
            maxcitenames=2,
            uniquename=init,
            isbn=false,
            backend=bibtex]{biblatex}
\usepackage{hyperref}

\newcommand{\field}{\mathbb}

\newcommand{\reals}{\field{R}}
\newcommand{\integers}{\field{Z}}

\newcommand{\R}{\reals}
\newcommand{\Z}{\integers}

\newcommand{\st}{\text{s.t.}}

\renewcommand{\th}{\text{th}}

\newcommand{\subprob}{\text{sub-problem}}

\newcommand{\set}[1]{\{#1\}}
\newcommand{\Set}[1]{\left\{#1\right\}}
\newcommand{\defset}[3][\defsep]{\set{#2#1#3}}
\newcommand{\Defset}[3][\defsep]{\Set{#2#1#3}}
\newcommand{\define}{\mathrel{{\mathop:}{=}}}

\makeatletter
\newcommand{\fcdot}{\,\cdot\,}
\newcommand{\fcarg}[1]{\def\fc@rg{#1}\ifx\fc@rg\empty\fcdot\else\fc@rg\fi}
\makeatother
\newcommand{\norm}[2][]{\lVert\fcarg{#2}\rVert\ifx#1\empty\else_{#1}\fi}
\newcommand{\Norm}[2][]{\left\lVert#2\right\rVert\ifx#1\empty\else_{#1}\fi}

\newcommand{\Abs}[1]{\left\lvert#1\right\rvert}

\newcommand{\upperbound}{U}
\newcommand{\ub}{\upperbound}
\newcommand{\lowerbound}{L}
\newcommand{\lb}{\lowerbound}

\DeclareMathOperator*{\argmax}{arg\,max}

\newcommand{\GammaValFun}[1]{\Phi_{\text{rob}}(#1)}
\newcommand{\SubProbValFun}[2]{\Phi_{#1}(#2)}
\newcommand{\robval}{v}

\newcommand{\subprobval}[1]{\robval_{#1}}

\newcommand{\feasupper}{X}

\newtheorem{assumption}{Assumption}

\definecolor{my-red}{HTML}{d7191c}

\newtheorem{lemma}{Lemma}
\newtheorem{theorem}{Theorem}
\newtheorem{proposition}{Proposition}

\bibliography{references}

\begin{document}

\title[Exact and Heuristic Methods for $\Gamma$-Robust Min-Max Problems]{Exact
  and Heuristic Methods for\\$\Gamma$-Robust Min-Max Problems}

\author[Y. Beck]{Yasmine Beck}

\address[Y. Beck]{%
  Eindhoven University of Technology,
  Department of Industrial Engineering and Innovation Sciences,
  PO Box 513,
  5600 MB Eindhoven,
  the Netherlands}%
\email{y.beck@tue.nl}

\date{\today}

\begin{abstract}
  Bilevel optimization is a powerful tool for modeling hierarchical
  decision-making processes, which arise in various real-world applications.
  Due to their nested structure, however, bilevel problems are
  intrinsically hard to solve---even if all variables are continuous and all
  parameters of the problem are exactly known. Further challenges arise if
  mixed-integer aspects and problems under uncertainty are considered.
  In this article, we summarize selected results from the author's dissertation
  \parencite{Beck:2024}. We study
  mixed-integer linear min-max problems with a $\Gamma$-robust treatment of
  uncertain data, for which we present exact and heuristic solution
  approaches. The performance of the methods is assessed in a computational
  study on~$560$~instances of the knapsack interdiction problem. Our results
  show that the heuristic closes the optimality gap for a
  significant portion of the considered instances and often practically
  outperforms both heuristic and exact benchmark approaches.
\end{abstract}

\keywords{Bilevel optimization, Robust optimization, Mixed-integer programming,
  Branch-and-Cut, Knapsack interdiction}
\subjclass[2020]{$90$C$11$, $90$C$27$, $90$C$57$, $90$C$70$}

\maketitle

\section{Introduction}
\label{sec:introduction}

Over the last years and decades, bilevel problems have gained increasing
attention because of their ability to model hierarchical interactions
between two decision-makers---the leader and the follower.
For an overview of the many applications of bilevel optimization, we refer
to \textcite{Dempe:2020} and to the recent surveys in \textcite{Kleinert_et_al:2021}
and \textcite{Beck_et_al:2023a}.
The latter focuses on bilevel problems under uncertainty, which is
also at the core of this article.
In what follows, we consider mixed-integer linear min-max problems of the form
\begin{subequations}
  \label{eq:min-max-prob}
  \begin{align}
    \min_{x,y} \quad
    & c^\top x + f^\top y
      \label{eq:min-max-prob:obj}
    \\
    \st \quad
    & x \in \feasupper,
      \label{eq:min-max-prob:ul}
    \\
    & y \in \argmax_{y'} \Defset{f^\top y'}{y' \in Y(x)}
      \label{eq:min-max-prob:ll}
  \end{align}
\end{subequations}
with~$Y(x) \subseteq \set{0,1}^{n_y}$,
$\feasupper \define \defset{x \in \R^{n_{\text{C}}} \times
  \Z^{n_{\text{D}}}}{Ax \geq a}$,
$n_x = n_{\text{C}} + n_{\text{D}}$, $c \in \R^{n_x}$, $f \in \R^{n_y}$,
$A \in \R^{l \times n_x}$, and~\mbox{$a \in \R^l$}.
We refer to~\eqref{eq:min-max-prob:obj}--\eqref{eq:min-max-prob:ul} as the
upper-level (or the leader's) problem and to~\eqref{eq:min-max-prob:ll} as the
lower-level (or the follower's) problem.
To ensure that an optimal solution to Problem~\eqref{eq:min-max-prob} exists,
we impose the following for the remainder of this article.
\begin{assumption}
  \begin{enumerate}
  \item For all~$x \in \feasupper$, the set~$Y(x)$ is non-empty.
  \item The set~$\defset{(x,y)}{x \in \feasupper,\, y \in Y(x)}$
    is non-empty and compact.
  \item All variables~$x$ that appear in the lower-level constraints are
    bounded integers.
  \end{enumerate}
\end{assumption}

In this article, we study Problem~\eqref{eq:min-max-prob} with uncertainty in
the objective function coefficients~$f$.
For all~$i \in [n_y] \define \set{1,\ldots,n_y}$, we thus
consider~$\bar{f}_i \in [f_i - \Delta f_i,f_i]$ instead of~$f_i$.
Here, $f_i$ is the nominal value of the~$i$\th\ objective function
coefficient and~\mbox{$\Delta f_i \geq 0$} is its maximum deviation from the
nominal value.
We address this kind of uncertainty using a~$\Gamma$-robust approach
\parencite{Bertsimas_Sim:2003} in which the follower hedges against at
most~\mbox{$\Gamma \in [n_y]$} deviations that adversely affect his optimal
objective function value.
This leads us to considering the bilevel problem
\begin{equation}
  \label{eq:rob-min-max}
  \min_{x,y} \quad c^\top x + f^\top y
  \quad \st \quad x \in X,\, y \in S_\Gamma(x),
\end{equation}
where~$S_\Gamma(x)$ is the set of optimal solutions to the
$\Gamma$-robust lower-level problem
\begin{equation}
  \label{eq:gamma-rc}
  \GammaValFun{x} \define
  \max_{y \in Y(x)} \Set{f^\top y -
  \max_{\defset{S \subseteq [n_y]}{|S|\leq \Gamma}}
  \sum_{i \in S}\Delta f_iy_i}.
\end{equation}
Using the optimal-value function~$\GammaValFun{x}$, we
can re-state Problem~\eqref{eq:rob-min-max} as
\begin{equation}
  \label{eq:rob-min-max-epi}
  \min_{x,\eta} \quad c^\top x + \eta
  \quad \st \quad x \in X,\, \eta \geq \GammaValFun{x}.
\end{equation}
In the dissertation \parencite{Beck:2024}, two solution approaches have been
derived---an exact branch-and-cut method and a heuristic---which are the first
to tackle Problem~\eqref{eq:rob-min-max} directly.
The methods have been published in \textcite{Beck_et_al:2023b} and
\textcite{Beck_et_al:2025}, respectively, and they rely on the
following auxiliary result. For further details and a proof of this result, we
refer to Lemma~1 and the respective discussion in \textcite{Beck_et_al:2025}.

\begin{lemma}
  \label{lem:ll-rc}
  Let~$x \in X$ be given arbitrarily and suppose that the indices are ordered
  such that~$\Delta f_i \geq \Delta f_{i+1}$ holds for all~$i \in [n_y]$
  with~$\Delta f_{n_y + 1} \define 0$.
  Then, the~$\Gamma$-robust counterpart~\eqref{eq:gamma-rc} of the lower-level
  problem can be solved by solving
  \begin{equation*}
    \GammaValFun{x} = \max_{\ell \in \mathcal{L}} \Set{\SubProbValFun{\ell}{x}},
  \end{equation*}
  where~$\mathcal{L} = \set{\Gamma + 1, \Gamma + 3, \Gamma + 5, \dots,
    \Gamma + \gamma, n_y + 1}$ with~$\gamma$ being the largest odd integer
  such that~$\Gamma + \gamma < n_y + 1$, and
  \begin{equation*}
    \SubProbValFun{\ell}{x} \define
    -\Gamma \Delta f_\ell + \max_{y \in Y(x)} \Set{%
      \sum_{i=1}^\ell (f_i - \Delta f_i + \Delta f_\ell)y_i
      + \sum_{i=\ell+1}^{n_y} f_iy_i}, \quad \ell \in \mathcal{L}.
  \end{equation*}
\end{lemma}

\section{An Exact Branch-and-Cut Approach}
\label{sec:bnc}

At the root node of the branch-and-cut search tree, we solve the linear problem
\begin{equation}
  \label{eq:master-prob}
  \min_{x,\eta} \quad c^\top x + \eta
  \quad \st \quad (x,\eta) \in \Omega_0
  \define \defset{(x',\eta') \in \bar{X} \times \R}{\eta' \geq \eta^-},
\end{equation}
which is obtained from Problem~\eqref{eq:rob-min-max-epi} by omitting the
constraint~$\eta \geq \GammaValFun{x}$ and by relaxing the integrality
restrictions for the leader's variables~$x$.
In~\eqref{eq:master-prob}, we use~$\bar{X}$ to denote the continuous relaxation
of~$X$.
Moreover, $\eta^- \in \R$ is a given lower bound on~$\GammaValFun{x}$ for
all~$x \in X$. Details on how to obtain such a bound are given in
\textcite{Beck_et_al:2023b}.
After considering Problem~\eqref{eq:master-prob},
we iteratively add valid inequalities or branch to separate
integer-infeasible points, and we also add valid inequalities to cut off
bilevel-infeasible points. At node~$k$ of the branch-and-cut search tree, we
consider the problem
\begin{equation}
  \label{eq:general-node-prob}
  \min_{x,\eta} \quad c^\top x + \eta
  \quad \st \quad (x,\eta) \in \Omega_k \subseteq \bar{X} \times \R.
\end{equation}
Here, $\Omega_k$ is obtained from~$\Omega_0$ by adding all valid
inequalities that have been generated at nodes along the path from the root
node to node~$k$ and by imposing all branching decisions that have been made
along that path.
If Problem~\eqref{eq:general-node-prob} is infeasible or if the
objective function value corresponding to an optimal solution~$(x^k,\eta^k)$
exceeds the current upper bound~$\ub$, we can prune node~$k$.
Otherwise, we proceed as follows.
First, we check if the leader's variables~$x^k$ satisfy all integrality
constraints ($x^k \in X$).
If this is not the case, we separate the current solution by either
exploiting standard cutting planes from mixed-integer programming
or by branching.
If~\mbox{$x^k \in X$} holds, we check for bilevel feasibility, that is
(i.e.), we check if~$\eta^k \geq \GammaValFun{x^k}$ is satisfied.
To this end, we need to solve the~$\Gamma$-robust lower-level
problem~\eqref{eq:gamma-rc}, which can be done by solving~$\Abs{\mathcal{L}}$
deterministic lower-level \subprob s; see Lemma~\ref{lem:ll-rc}.
If there is at least one such \subprob~$\ell \in \mathcal{L}$ for
which~$\eta^k < \SubProbValFun{\ell}{x^k}$ holds, the current
solution~$(x^k,\eta^k)$ is not bilevel feasible and we add a cut to
separate this point.
Problem-tailored cuts for the important class of monotone
interdiction problems that can be used for this purpose have been derived in
the author's dissertation; see Section~3.3 in~\textcite{Beck:2024} for the
details.
To sum up, the method to process node~$k$ of the branch-and-cut search tree
is formally stated in Algorithm~\ref{alg:node-processing}.
\begin{algorithm}
  \begin{algorithmic}[1]
    \STATE Solve Problem~\eqref{eq:general-node-prob}.
    \label{alg:step1}
    \IF{Problem~\eqref{eq:general-node-prob} is infeasible}
    \STATE Prune the current node and return to the main method.
    \ENDIF
    \STATE Let~$(x^k,\eta^k)$ denote the solution to
    Problem~\eqref{eq:general-node-prob}.
    \IF{$c^\top x^k + \eta^k \geq \ub$}
    \STATE Prune the current node and return to the main method.
    \ENDIF
    \IF{$x^k \notin X$}
    \STATE Either generate a cut valid for
    $\Omega_k \cap (X \times \R)$, augment $\Omega_k$, and go to
    Line~\ref{alg:step1}, or branch.
    \ENDIF
    \FORALL{$\ell \in \mathcal{L}$}
    \STATE \label{alg:generic-subprob}
    Compute~$\SubProbValFun{\ell}{x^k}$.
    \IF{$\eta^k < \SubProbValFun{\ell}{x^k}$}
    \STATE \label{alg:generic-cut}
    Generate a valid cut that excludes~$(x^k,\eta^k)$ from $\Omega_k$ and
    augment $\Omega_k$.
    \ENDIF
    \ENDFOR
    \STATE Set $\GammaValFun{x^k} \gets
    \max_{\ell \in \mathcal{L}} \set{\SubProbValFun{\ell}{x^k}}$ and $\ub \gets
    \min \set{U,\, c^\top x^k + \GammaValFun{x^k}}$.
    \STATE If at least one cut has been added in Line~\ref{alg:generic-cut},
    go to Line~\ref{alg:step1}.
  \end{algorithmic}
  \caption{Processing Node~$k$ of the Branch-and-Cut Search Tree}
  \label{alg:node-processing}
\end{algorithm}

\begin{theorem}[See Theorem~1 in \textcite{Beck_et_al:2023b}]
  If we embed Algorithm~\ref{alg:node-processing} in a classic
  branch-and-bound framework, we obtain a method that terminates with a
  globally optimal solution~$(x^*,\eta^*)$ to
  Problem~\eqref{eq:rob-min-max-epi} after investigating a finite number
  of nodes and after adding an overall finite number of cuts.
\end{theorem}

\section{A Heuristic in the Spirit of Bertsimas \& Sim}
\label{sec:heuristic}

Given the overall hardness of~$\Gamma$-robust min-max problems, which are
$\Sigma_2^p$-hard in general, we also derive a heuristic for these problems.
The method relies on the solution of a linear number of appropriately chosen
deterministic min-max problems. More formally, we have the following result.
For a proof of this result, we refer to Proposition~1
in~\textcite{Beck_et_al:2025}.

\begin{proposition}
  For all~$\ell \in \mathcal{L}$, let
  $\subprobval{\ell} \define \min_{x \in X} \set{c^\top x + \SubProbValFun{\ell}{x}}$.
  Then, $\subprobval{\ell}$ is a valid lower bound for
  the optimal objective function value of Problem~\eqref{eq:rob-min-max}.
\end{proposition}

The heuristic for Problem~\eqref{eq:rob-min-max} is formally stated in
Algorithm~\ref{alg:heuristic}.
The method starts by solving~$\Abs{\mathcal{L}}$ deterministic min-max
problems. Afterward, we use the solutions~$(x^\ell)_{\ell \in \mathcal{L}}$ to
these problems to compute upper bounds for Problem~\eqref{eq:rob-min-max}.
\begin{algorithm}
  \begin{algorithmic}[1]
    \STATE Set~$x^* \gets \textsf{None}$, $\lb \gets - \infty$,
    and~$\ub \gets \infty$.
    \FORALL{$\ell \in \mathcal{L}$}
    \STATE Compute a solution~$x^{\ell}$ to the deterministic min-max problem
    \vspace*{-0.5em}
    \begin{equation*}
      \subprobval{\ell} \gets \min_{x \in X} \Set{c^\top x +
        \SubProbValFun{\ell}{x}}.
    \end{equation*}
    \vspace*{-1.5em}
    \label{alg:deterministic-min-max-subprob}
    \ENDFOR
    \STATE Set~$\lb \gets \max_{\ell \in \mathcal{L}}
    \set{\subprobval{\ell}}$ and~$i \gets 1$.
    \WHILE{$i \leq \Abs{\mathcal{L}}$ and~$\lb < \ub$}
    \label{alg:begin-while}
    \STATE Use Lemma~\ref{lem:ll-rc} to compute~$\GammaValFun{x^{\ell_i}}$.
    \IF{$c^\top x^{\ell_i} + \GammaValFun{x^{\ell_i}} < \ub$}
    \STATE Set~$x^* \gets x^{\ell_i}$ and~$\ub \gets c^\top x^* +
    \GammaValFun{x^*}$.
    \ENDIF
    \STATE Set~$i \gets i + 1$.
    \label{alg:end-while}
    \ENDWHILE
    \RETURN $x^*,\, \lb,\, \ub$
  \end{algorithmic}
  \caption{Heuristic for~$\Gamma$-Robust Min-Max Problems}
  \label{alg:heuristic}
\end{algorithm}

\begin{theorem}[See Theorem~1 in~\textcite{Beck_et_al:2025}]
  Algorithm~\ref{alg:heuristic} returns a feasible leader's decision~$x^*$ as
  well as valid lower and upper bounds~$\lb$ and~$\ub$
  for Problem~\eqref{eq:rob-min-max}.
\end{theorem}

In \textcite{Beck:2024,Beck_et_al:2025}, we further derive sufficient
conditions under which Algorithm~\ref{alg:heuristic} terminates with a provably
globally optimal solution after Line~\ref{alg:deterministic-min-max-subprob},
i.e., after only solving deterministic min-max problems.
Note that this extends the famous result by Bertsimas and Sim
\parencite{Bertsimas_Sim:2003} to the~$\Gamma$-robust min-max setting.

\section{Computational Results}
\label{sec:computational-results}

\begin{figure}
  \begin{center}
    \input{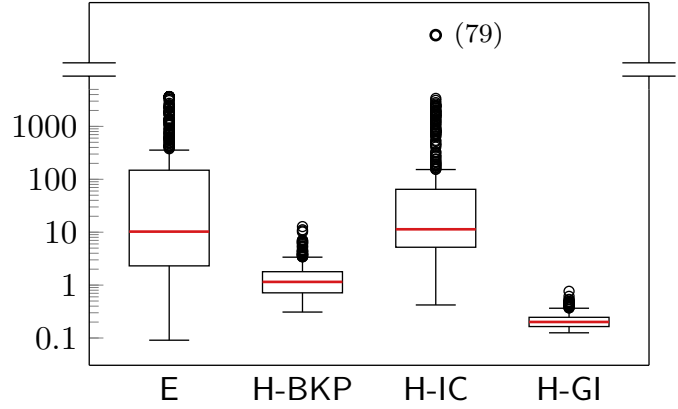}
    \caption{Box-plots of the runtimes for~\textsf{E}, \textsf{H-BKP},
      \mbox{\textsf{H-IC}}, and~\textsf{H-GI}. Runtimes
      (in~\si{\second}) are depicted on a log-scaled~$y$-axis.}
    \label{fig}
  \end{center}
\end{figure}

In this article, we report numerical results for~$560$ instances of
the~$\Gamma$-robust knapsack interdiction problem with continuous
deviations~$\Delta f$. Details regarding the generation of the instances and
the computational setup can be found in \textcite{Beck_et_al:2025}.
We compare four solution approaches.
The first is the exact branch-and-cut method presented in
Section~\ref{sec:bnc} in which we use the problem-tailored cuts
derived in \textcite{Beck_et_al:2023b}. We refer to this method as \textsf{E}.
Moreover, we consider two variants of the heuristic in
Algorithm~\ref{alg:heuristic}---one using the \textsf{bkpsolver}
\parencite{Weninger_Fukasawa:2022} and one using the branch-and-cut method in
\textcite{Fischetti_et_al:2019} to solve the deterministic min-max problems.
We refer to these approaches as \textsf{H-BKP} and \textsf{H-IC},
respectively.
Finally, we compare our methods to the ``Greedy Interdiction''
heuristic presented in \textcite{DeNegre:2011}, which we abbreviate as
\textsf{H-GI}.
Table~\ref{tab} and Figure~\ref{fig} summarize our numerical results.

\begin{table}
  \centering
  \caption{The number of instances for which a feasible point with finite
      gap is found (``feasible'') and the number of instances solved to global
      optimality (``optimal'') for~\textsf{E}, \textsf{H-BKP}, \textsf{H-IC},
      and~\textsf{H-GI}.
      For those instances with finite but non-zero gap (``open gap''), also the
      average gap (``average gap''; in \si{\percent}) is shown.
    }
    \begin{tabular}{lrrrr}
      \toprule
      & feasible & optimal & open gap & average gap\\
      \midrule
      \textsf{E}     & 560 & 524 &  36 &   7.03\\
      \textsf{H-BKP} & 560 & 554 &   6 &   0.08\\
      \textsf{H-IC}  & 481 & 476 &   5 &   0.10\\
      \textsf{H-GI}  & 560 &   4 & 556 & 100.00\\
      \bottomrule
    \end{tabular}
  \label{tab}
\end{table}

All methods except for \textsf{H-IC} find feasible points for all~$560$
instances within~\SI{1}{\hour}. On the subset of instances
that \textsf{H-IC} can tackle, \textsf{H-IC} performs slightly better than
\textsf{E} in terms of runtimes.
Overall, \textsf{H-GI} achieves the smallest runtimes, but its solution quality
is rather poor.
In contrast, \textsf{H-BKP} not only outperforms \textsf{E} by significant
orders of magnitude in terms of runtimes, it also proves global optimality for
almost all considered instances.
These results indicate that, if efficient black-box solvers are available for
the deterministic min-max problems, our heuristic can outperform both exact
and heuristic benchmark approaches.

\section{Summary}
\label{sec:summary}

In this article, we summarize selected results from the author's
dissertation~\parencite{Beck:2024}. To this end, we present exact and heuristic
solution approaches for mixed-integer linear min-max problems with
a~$\Gamma$-robust treatment of objective uncertainty.
The performance of the methods is assessed in a computational study
on~$560$~instances of the knapsack interdiction problem. Our results show that
the heuristic closes the optimality gap for a significant portion of the
considered instances and often practically outperforms both heuristic and exact
benchmark approaches.


\printbibliography

\end{document}